\documentclass{article}

\usepackage{graphicx}      
\usepackage{amsmath}
\usepackage{amssymb}
\newcommand{\ltwo}{{\ell^2}}

\title{Modal approach to the controllability problem of distributed parameter systems with damping}

\author{Alexander Zuyev}
\date{\small 
 Institute of Applied Mathematics and Mechanics, \\ National Academy of Sciences of Ukraine}

\begin{document}
\maketitle

\begin{abstract}                
This paper is devoted to the controllability analysis of a class of linear control systems in a Hilbert space.
It is proposed to use the minimum energy controls of a reduced lumped parameter system for solving the infinite dimensional steering problem approximately.
Sufficient conditions of the approximate controllability are formulated for a modal representation of a flexible structure with small damping.\footnote{This work was supported by grants from the National Academy of
Sciences of Ukraine and the Alexander von Humboldt Foundation.}
\end{abstract}



\section{Introduction}
The problems of spectral, approximate, exact, and null controllability of distributed parameter systems have been intensively studied over the last few decades~\cite{Balakrishnan,Coron,LL}. On the one hand, the question of the approximate controllability of a linear time-invariant system on a Hilbert space can be formulated in terms of an invariant subspace of the corresponding adjoint semigroup~\cite{LR},~\cite[p.~56]{Coron}.
On the other hand, the problem of effective control design remains challenging for wide classes of mechanical systems (see, e.g.~\cite{Z2003,ZuCDC2006,Z2006,ZuMPE,Z2016}, and references therein).
The goal of this work is to propose a constructive control strategy, based on a reduced model, and to justify that this approach can be used to solve the approximate controllability problem for infinite dimensional systems.

\section{Problem statement}
This paper addresses the problem of approximate controllability of a linear differential equation
\begin{equation}
\dot x = Ax + Bu,\quad x\in H,\;u \in {\mathbb R}^m,
\label{l2_eq}
\end{equation}
where $H$ is a Hilbert space, $A:D(A)\to H$ is a closed densely defined operator, and $B:{\mathbb R}^m\to H$ is a continuous operator. We assume that $A$ generates a strongly continuous semigroup of operators $\{e^{tA}\}_{t\ge 0}$ on $H$. Hence, for any $x^0\in H$ and $u\in L^2(0,\tau)$, the mild solution of~\eqref{l2_eq} corresponding to the initial condition $x|_{t=0}=x^0$ and control $u=u(t)$ can be written as follows:
\begin{equation}
x(t;x^{0},u) = e^{tA}x^{0} + \int_0^t e^{(t-s)A}Bu(s)\,ds,\; 0\le
t\le \tau.
 \label{mildsol}
\end{equation}

Note that that system~\eqref{l2_eq} is approximately controllable in time $\tau>0$ if~(cf.~\cite{Coron}), given $x^{0},x^{1}\in H$ and $\varepsilon>0$,
there exists $u\in L^2(0,\tau)$ such that
$
\|x(\tau;x^{0},u)-x^{1}\|<\varepsilon
$.
In order to study the approximate controllability of system~\eqref{l2_eq}, we use the following result.

{\bf Proposition~1.}~\cite{Zu2009}
{\em Let $\{Q_N\}_{N=1}^\infty$ be a family of bounded linear
operators on $H$ satisfying the following conditions:
\begin{enumerate}
\item[1)]
\begin{equation}
\lim_{N\to\infty}\|Q_N x\|= 0,\quad \text{for all} \; x\in H;
\label{Q_cond}
\end{equation}
\item[2)] the operators $e^{tA}$ and $Q_N$ commute;
\item[3)] for each $x^{0},x^{1}\in H$, $N\ge 1$, there is a control
$u^{N}_{x^{0},x^{1}}\in L^\infty(0,\tau)$ such that
\begin{equation}
(I-Q_N)\Bigl(x(\tau;x^{0},u^{N}_{x^{0},x^{1}})-x^{1}\Bigr) =0,
\label{proj_cond}
\end{equation}
\begin{equation}
\lim_{N\to\infty} \left(\|Q_N B\|\cdot
\|u^{N}_{x^{0},x^{1}}\|_{L^2(0,\tau)}\right) =0. \label{BN_cond}
\end{equation}
\end{enumerate}
Then system \eqref{l2_eq} is approximately controllable in
time $\tau$, and the above family of functions
$u=u^{N}_{x^{0},x^{1}}(t)$, $0\le t\le \tau$, can be used to solve
the approximate controllability problem. }

For a possible application of this proposition, we assume that each operator $P_N=I-Q_N$ is a finite dimensional projection. Let ${\rm dim }({\rm Im} \,P_N) = d_N$. For given $x^{0},x^{1}\in H$, we introduce vectors
$$
\tilde x_N^0=P_N x^0,\;\tilde x_N^1=P_N x^1,\;\tilde x_N = P_N x,\quad (\tilde x_N^0,\tilde x_N^1,\tilde x_N\in {\rm Im} \,P_N),
$$
and operators $\tilde A_N = P_N A$, $\tilde B_N = P_N B$. Then condition~\eqref{proj_cond} implies that $u^{N}_{x^{0},x^{1}}(t)$ should solve the following control problem:
\begin{equation}
\dot{\tilde x}_N = \tilde A_N \tilde x_N + \tilde B_N u,\quad t\in [0,
\tau],
\label{ODE_approx}
\end{equation}
$$
\tilde x_N|_{t=0}=\tilde x_N^0,\;\tilde x_N|_{t=
\tau}=\tilde x_N^1.
$$
Here we have used the assumption that $P_N$ and $A$ commute as well as the property $P_N=P_N^2$ of a projection. To satisfy condition~\eqref{BN_cond}, it is natural to look for a control $u=u^{N}_{x^{0},x^{1}}(t)$ that minimizes the functional
\begin{equation}
J = \int_0^\tau (Qu,u)\,dt \to \min
\label{cost}
\end{equation}
with some symmetric positive definite $m\times m$-matrix $Q$.
As control system~\eqref{ODE_approx} evolves on a real $d_N$-dimensional vector space ${\rm Im}\, P_N$, we may treat~\eqref{ODE_approx} as a system on ${\mathbb R}^{d_N}$ without lack of generality. By applying the Pontryagin maximum principle, we get the optimal control for problem~\eqref{ODE_approx}-\eqref{cost}:
\begin{equation}\
\tilde u(t) = Q^{-1}\tilde B_N'e^{(\tau-t)\tilde A_N'}\nu,\quad
\nu = \left(\int_0^\tau e^{s \tilde A_N}\tilde B_N Q^{-1}\tilde B_N' e^{s \tilde A_N'}\,ds\right)^{-1} (\tilde x^1_N-e^{\tau \tilde A_N}\tilde x^0_N),
\label{optimalcontrol}
\end{equation}
where the prime stands for the transpose.
Proposition~1 implies that the proof of the approximate controllability
can be reduced to the checking conditions~\eqref{Q_cond} and~\eqref{BN_cond} with a family of smooth controls $u^{N}_{x^{0},x^{1}}=\tilde u(t)$ given by~\eqref{optimalcontrol}.
The main contribution of this paper is the application of such a scheme for a class of systems~\eqref{l2_eq} representing the oscillations of a flexible structure with damping.

\section{Flexible system with damping}
Consider a particular case of system~\eqref{l2_eq}
as follows
\begin{equation}
\dot x = Ax + Bu,\quad x=(x_1,x_2,...)'\in \ltwo,\;u\in {\mathbb R},
\label{l2e}
\end{equation}
where
$
\|x\|_\ltwo = \left(\sum_{n=1}^\infty x_n^2\right)^{1/2}
$. We assume that the operator $A:D(A)\to\ltwo$ in~\eqref{l2e} is given by its block-diagonal matrix:
$$
A={\rm diag}(A_0,A_1,A_2,...),\;
A_0 = \left(
\begin{array}{cc}
0 & 1 \\
0 & 0
\end{array}
\right),\;
A_n = \left(
\begin{array}{cc}
0 & \omega_n \\
-\omega_n & -2\kappa
\end{array}
\right),\quad n=1,2,... \;,
$$
and
$
B=(0,1,0,b_1,0,b_2,...)'\in\ltwo
$.
Control system~\eqref{l2e} is a linear model of a rotating flexible beam attached
to a rigid body. The components of $x$ plays the role of modal coordinates, and the control $u$ is the angular acceleration of the body. Coefficients $\omega_n$ and $b_n$ are, respectively, the modal frequency and the control coefficient corresponding to the $n$-th mode of oscillations of the beam. The coefficient $\kappa>0$ represents the viscous damping in the beam.
The procedure of deriving the equations of motion with modal coordinates is described in the paper~\cite{ZuAutom} for a rotating rigid body with flexible beams.

The main result of this paper is as follows.

{\bf Proposition~2.}{\em
Assume that $b_n\neq 0$ and $\omega_n>0$ for all $n=1,2,...$.
Then there exists a $\tau>0$ such that system~\eqref{l2e}
is approximately controllable in time $\tau$ provided that
\begin{equation}
\sum_{\scriptsize\begin{array}{c}i,j=1\\i\neq j\end{array}}^\infty
\frac{1}{(\omega_i-\omega_j)^2}<\infty \label{seriesconv}
\end{equation}
and that the damping coefficient $\kappa$ is small enough.
}

{\bf Proof.}
Let us introduce the family of operators $P_N:
\ltwo\to\ltwo$ as follows:
$$
P_N\left(\begin{array}{c}
\xi_0 \\ \eta_0 \\  \vdots \\ \xi_N \\ \eta_N \\ \xi_{N+1} \\ \eta_{N+1} \\ \vdots
\end{array}\right) = \left(\begin{array}{c}
\xi_0 \\ \eta_0 \\  \vdots \\ \xi_N \\ \eta_N \\ 0 \\ 0 \\ \vdots
\end{array}\right),
$$
and $Q_N=I-P_N$, $N=1,2,...$. Then condition~1) of Proposition~1 holds. To check condition~2), we compute the semigroup $\{e^{tA}\}$ generated by $A$:
$$
e^{tA}={\rm diag}\left(e^{tA_0},e^{tA_1},e^{tA_2},...\right),
$$
$$
e^{tA_0}=\left(\begin{array}{cc}1 & t \\ 0 & 1\end{array}\right),\quad
e^{tA_n}=e^{-\kappa t}\left(\begin{array}{cc}1 & t \\ 0 & 1\end{array}\right),
$$

The assertion of Proposition~2 follows from
Proposition~1 by exploiting the construction of $L^2$-minimal
controls $u^{N}_{x^{0},x^{1}}=\tilde u(t)$
in~\eqref{optimalcontrol}.

\section{Conclusions}
This work extends the result of~\cite{Zu2009} for the case of a flexible system with damping. As it was shown earlier in~\cite{Zu2009}, condition~\eqref{seriesconv} is satisfied for the Euler-Bernoulli beam without damping. Hence, condition~\eqref{seriesconv} is sufficient
for the approximate controllability in both conservative ($\kappa=0$) and dissipative (small $\kappa>0$) cases under our assumptions.
An open question is whether is it possible to relax restrictions~\eqref{seriesconv} in order to justify the relevance of
controls~\eqref{optimalcontrol} under a weaker assumption on the distribution of  the modal frequencies $\{\omega_n\}$.

\end{document}